\documentclass[nobib]{tufte-handout}
\usepackage{amsmath,stmaryrd,amssymb,amsthm,url,booktabs,hyperref}

\title{Explicit String bundles\thanks{\href{https://orcid.org/0000-0002-3478-0522}{orcid.org/0000-0002-3478-0522}\\
\noindent This document is released under a CC0 license: \href{http://creativecommons.org/publicdomain/zero/1.0/}{\texttt{creativecommons.org/publicdomain/zero/1.0/}}.}}
\author{David Michael Roberts}
\date{1 July 2014}

\usepackage{euler}

\usepackage[style=alphabetic,%
            citestyle=alphabetic]{biblatex}
\renewcommand{\a}{\alpha}
\renewcommand{\b}{\beta}
\newcommand{\g}{\gamma}
\renewcommand{\d}{\delta}
\newcommand{\e}{\varepsilon}
\renewcommand{\l}{\lambda}
\newcommand{\HP}{\mathbb{HP}^1}
\newcommand{\HH}{\mathbb{H}}
\newcommand{\sourcetarget}{\rightrightarrows}
\DeclareMathOperator{\String}{String}
\DeclareMathOperator{\Spin}{Spin}

\theoremstyle{definition}
\newtheorem*{challenge}{Challenge}
\newtheorem*{note}{Note}
\newtheorem*{proposition}{Proposition}
\newtheorem*{exercise}{Exercise}

\begin{document}
\maketitle

\marginnote{These notes are the written version of a talk delivered at Herriot-Watt University on 26 June 2014 at the Workshop on Higher Gauge Theory and Higher Quantization. DMR is supported by ARC grant number DP120100106.\\
Thanks to David Baraglia, Michael Murray, Christian Saemann and Raymond Vozzo for helpful converstations.}
If we look back at the historical development of bundles, the notion
of a principal $H$-bundle for $H$ a Lie group arose via considerations
of homogeneous spaces $G/H$, and the defining bundle $G\to G/H$. Here
$H$ is a closed subgroup of $G$, and this will be a general assumption
for the talk. For instance, we can conside Stiefel manifolds, Grassmann
manifolds, projective spaces, Minkowski space, spheres,\ldots. I am going
to leverage this in order to address the

\begin{challenge}
Write down a (nontrivial) 2-bundle. Equivalently, 
write down a \v Cech cocycle with values in an interesting crossed 
module $(K\xrightarrow{t}H,H\times K \xrightarrow{a} K)$. 
\end{challenge}

Recall \cite{Breen94} that the cocycle equations are
\begin{align*}
	h_{ij}^\a h_{jk}^\b =& t(k_{ijk}^{\a\b\g}) h_{ik}^\g  \\ 
	t(h_{ij}^a,k_{jkl}^{\b\l\e}) k_{ijl}^{\a\e\d}  =& k_{ijk}^{\a\b\g} k_{ikl}^{\g\l\d}.
\end{align*}
where the $h_{ij}^\a$ are $H$-valued functions, the $k_{ijk}^{\a\b\g}$ 
are $K$-valued functions and the two sorts of indices label open sets 
of the base space. We shall return to this momentarily. Note that at this
point we haven't even started to consider connections, which are necessary
for gauge theory (and in fact we won't even go so far today).

\begin{note}
I am \emph{not} going to use good open covers (that is,
those such that non-empty finite intersections are contractible), since
in many geometric situations there are naturally arising open covers
that are not good. Instead, I will be using \emph{truncated globular
hypercovers} (these are open covers with particular properties), and I will
define these in a moment. For now, suffice it to say, this is why there are
two different sorts of indices on the cocycle.
\end{note}

Christian Saemann asked (Feb 2013): I want a 2-bundle on (conformally 
compactified)\footnote{recall that this is diffeomorphic to $S^5\times S^1$} $\mathbb{R}^{5,1}$.
So let's try lifting the frame bundle of $S^5\times S^1$ to a \emph{String}
bundle. Note that the $S^1$ factor contributes nothing (its frame bundle is 
trivial) so just work over $S^5$. Note that the frame bundle of $S^5$
is most definitely not trivial.

The frame bundle $FS^5 \to S^5$ is classified by a map $S^5 \supset S^4
\to SO(5)$, called the \emph{clutching} or \emph{transition function}.
Since $S^5$ is 4-connected, the first Stiefel-Whitney class $w_1$ necessarily
vanishes, as does the characteristic class $p_1/2 \in H^4(S^5,\mathbb{Z})$
that is the obstruction to lifting to a String bundle. Thus we can be 
assured that the lift we are after does exist. From the vanishing
of $w_1$ we know the transition function lifts to a function $S^4 \to \Spin(5)$,
and so defines a class in $\pi_4(\Spin(5))$, which is the group $\mathbb{Z}/2\mathbb{Z}$
\cite{Mimura-Toda_64}. Since $FS^5$ is not trivial, the transition function
needs to represent the non-trivial homotopy class. We want to write 
down an explicit function in coordinates, rather than use some abstract 
representative.

To approach this, we first use the exceptional isomorphism $\Spin(5) \simeq Sp(2)$,
where $Sp(2)$ is the group of $2\times2$ unitary quaternionic matrices.
The non-trivial class in $\pi_4(Sp(2))$ is represented by a map
$S^4 \to S^3 \simeq Sp(1) \hookrightarrow Sp(2)$ and here $Sp(1)$ is the group
of unit quaternions. The map between spheres is (up to homotopy) the 
suspension of the Hopf map $S^3 \to S^2$, which is not a priori a smooth map, 
and the inclusion is $q \mapsto \left(\begin{array}{cc}q&0\\0&1 \end{array}\right)$.
Note that this implies that $FS^5$ lifts to an $Sp(1)$-bundle, and this is
what we shall assue without further comment.

The first task is then to write down a smooth, non-null-homotopic smooth map
$S^4 \to Sp(1)$. We shall use quaternionic coordinates on $S^4 = \HP$, that is,
homogeneous coordinates $[p;q]$ where at least one of $p$, $q$ is non-zero.

\begin{proposition}
The smooth function
\[
	T[p;q] = \frac{2p\bar q i \bar p q - |p|^4 + |q|^4}{|p|^4 + |q|^4}\quad (\in Sp(1))
\]
represents the non-trivial class of $\pi_4(Sp(1))$, and hence is the transition 
function for $FS^5$.\end{proposition}

Now we want to shift perspective a little bit, and note that the function
$T$ gives rise to a smooth \emph{functor} from the \v Cech groupoid $U\times_{S^5} U 
\sourcetarget U$ over $S^5$ coming from the open cover by two discs\footnote{One should take these as open discs,
and so the intersection would be $S^4 \times (-\e,\e)$; we extend $T$ to this 
slightly larger subspace by taking it constant in the direction of the interval.} 
$U := D_+ \coprod D_- \to S^5$.
For future notational convenience, write $U^{[2]} = U\times_{S^5} U$.

Since we now have an explicit \v Cech cocycle (this is precisely what the above
functor is) for $FS^5$, we can talk about lifting this to a \v Cech cocycle for
the 2-group $\String_{Sp(1)}$. But what \emph{is} this? There are many models for
String 2-groups, and we shall take the crossed module $(\widehat{\Omega Sp(1)} \to PSp(1))$,
where $PSp(1)$ is the group of smooth paths $[0,1] \to Sp(1)$ based at $1\in Sp(1)$,
and $\widehat{\Omega Sp(1)}$ is the universal central extension of the subgroup
$\Omega Sp(1) \subset PSp(1)$ of loops \cite{BCSS_07}.
Note that the abstract details of what I'm considering doesn't rely on this
choice of model. Notice that $(\widehat{\Omega Sp(1)} \to PSp(1))$ comes with a map to the crossed
module $(1\to Sp(1))$, and that the former gives rise to a groupoid (which I
shall call $\String(3)$, as $Sp(1) \simeq \Spin(3)$), namely
the action groupoid for $\widehat{\Omega Sp(1)}$ acting on $PSp(1)$ via the
given homomorphism, and a 2-groupoid $\mathbf{B}\String(3)$ with a single object 
(using the 2-group structure). More generally, we can repeat these constructions
with any compact, simple, simply connected Lie group $G$ to get a 2-group $\String_G$.
Also, given an inclusion of Lie groups\footnote{that induces an isomorphism
$H^3(G,\mathbb{Z}) \to H^3(H,\mathbb{Z})$; the examples listed below all satisfy this, as can be calculated via the long-exact sequence in homotopy.} $H\to G$ gives an inclusion of Lie 2-groups 
$\String_H \to \String_G$.

In the \v Cech groupoid $U^{[2]} \sourcetarget U$ we don't have $U^{[2]}$
a disjoint union of contractible opens, so we take an open cover $V\to U^{[2]}$
where $V$ \emph{is} such a disjoint union (or, at least, acyclic enough). 
Since the non-trivial part of $U^{[2]}$
is $D_+\cap D_- \sim \HP$, we will take $V$ to be the two $\HH$ charts
$\HH_+$ and $\HH_-$ given by non-vanishing of each of the two homogenous coordinates.
Then if we take the fibred product $V^{[2]} = V\times_{U^{[2]}}V$ we get a Lie 2-groupoid
$V^{[2]} \sourcetarget V \sourcetarget U$, which I call a \emph{truncated globular 
hypercover}.\footnote{This may look familiar if you've seen bundle 2-gerbes \cite{Stevenson_04} before.} The
nontrivial component of $V^{[2]}$ (it contains boring bits like $D_+$) is
the intersection $\HH_+\cap \HH_- = \HH^\times$. 
Notice that if we wanted to use a good open cover then $U$ would necessarily 
have had more open sets, and so more overlaps. In some sense we have made a 
trade-off in the number of open sets and the slight increase in complexity of the 
description. Also, we can finally see where the two sorts of indices in the
cocycle equation above come from: the indices $i,j,\ldots$ label open sets appearing
in $U$, and the indices $\a,\b,\ldots$ label the open sets appearing in $V$.

So, finally, a \v Cech cocycle on $S^5$ with values in $\String(3)$ is `just' a
2-functor
\[
	(V^{[2]} \sourcetarget V \sourcetarget U) \to \mathbf{B}\String(3).
\]
If we break this down, it is determined by components
\begin{align*}
	& V \to PSp(1)\\
	& V^{[2]} \to \widehat{\Omega PSp(1)} 
\end{align*}
and since we have so few open sets in the globular hypercover, functoriality
follows automatically. In our particular case, we want the first map to lift
the given $V \to U^{[2]} \to Sp(1)$.

Recall that $V$ is (essentially) $\HH_+ \coprod \HH_-$, we define the lift
in two parts:
\begin{align*}
	T_+(q) & = \left( s \mapsto \frac{|q|^4 - s^2 + 2\bar q i q}{|q|^4 + s^2} \right)\\
	T_-(p) & = \left( s \mapsto \frac{|p|^4s^2 - 1 + 2\bar p i p}{|p|^4s^2 + 1}\cdot \left(\frac{s-i}{s+i}\right)^2  \right)
\end{align*}
To define the remaining component of the 2-functor, we first take the difference of these two
maps to get a function $\HH^\times \to \Omega Sp(1)$
\[
	T_\Omega(q) = \left(s\mapsto \frac{(s + Q)(sQ - 1)}{(s - Q)(sQ + 1)} \cdot \left(\frac{s - i}{s + i}\right)^2 \right),
		\quad\text{where } Q = \bar q i q.
\]
Now we need to lift this map through the projection $\widehat{\Omega Sp(1)} 
\to \Omega Sp(1)$ (this is not a priori possible, but one calculates
the possible obstructions and they vanish). To do this, we need a workable
description of what $\widehat{\Omega Sp(1)}$ is. There are multiple papers
constructing this e.g.\ \cite{Mickelsson_87,Murray88,Mur_Stev03}. We
shall use the description of it as the quotient group
\[
	\frac{P\Omega Sp(1) \rtimes U(1)}{\widetilde{\Omega^2Sp(1)}}
\]
The precise embedding of the simply-connected covering group $\widetilde{\Omega^2Sp(1)}$
is not important, just that we can represent elements as equivalence
classes of pairs consisting of paths in $\Omega Sp(1)$ and elements of $U(1)$.

One calculates the final answer to be as follows. For any $q\in \HH^\times$, let
$q_t$ be any path (in $\HH^\times$) $1 \rightsquigarrow q$, and the lift to the central
extension is
\[
	T_{\widehat{\Omega}}(q) = [T_\Omega(q_t),1].
\]
This is independent of the choice of path and is smooth. This function, together
with $T_{\pm}$, defines the \v Cech cocycle we are interested in. We know that
this cocycle is not a coboundary, since geometrically realising everything we get a 
map $S^5 \to B\String(3)$ that picks out the nontrivial class in $\pi_5(B\String(3))
\simeq \pi_5(B\Spin(3)) \simeq \pi_4(\Spin(3)) \simeq \pi_4(Sp(1)) = \mathbb{Z}/2\mathbb{Z}$. One can also check (easily, as there are so few open sets involved
in the open covers), that these functions satisfy the cocycle equations 
displayed at the beginning of the notes.

Now this is just one example, and a pretty exceptional example at that, as the
dimensions involved are right on the boundary of where the obstructions vanish, not
to mention the use of quaternions. One can take a more global approach that
leads to many more examples as follows. The total space of the frame bundle $FS^5$, 
as an $Sp(1)$-bundle, is nothing other than the homogenous bundle $SU(3) \to SU(3)/Sp(1)=S^5$,
using the embedding $Sp(1) \simeq SU(2) \to SU(3)$ as a block matrix. One can calculate that 
$\String_{SU(3)}/\String(3) \simeq SU(3)/Sp(1)$, so that the underlying groupoid
of $\String_{SU(3)}$ is the `total space' of the $\String(3)$ bundle. Another way 
to view this is to consider the transitive $\String_{SU(3)}$ action on $S^5$ via 
the projection to $SU(3)$; then $\String(3)$ is the stabiliser of the basepoint. 

This picture generalises to any $\String_G$ acting on $G/H$ for $H < G$, and at this
point we can use any model of $\String_G$, including non-strict models, and even
2-groups in differentiable stacks, which have underlying Lie groupoids. There 
are a number of interesting exceptional examples which should be amenable to
the same treatment as above, for instance:
\begin{itemize}
	\item $\String_{G_2} \to G_2/SU(3) = S^6$
	\item $\String_{\Spin(7)} \to \Spin(7)/G_2 = S^7$
	\item $\String_{Sp(2)} \to Sp(2)/Sp(1) = S^7$
	\item $\String_{F_4} \to F_4/\Spin(9) = \mathbb{OP}^2$
\end{itemize}
The first three of these have explicit transition functions written down by
P\"uttmann in \cite{Puettmann_11}. $\mathbb{OP}^2$ admits a cover by three 
$\mathbb{R}^{16}$ charts, and is 7-connected. 

\begin{exercise}
Write down
transition functions for the $\Spin(9)$ bundle on $\mathbb{OP}^2$, and lift
them to $\String(9) = \String_{\Spin(9)}$-valued transition functions using a 
globular hypercover.
\end{exercise}

The astute reader will have realised that this method only gives a single
example on each homogeneous space with that particular structure group, which
in the case of $S^5$ is ok as there is only one nontrivial $\String(3)$ bundle.
But, for instance, $\String_{SU(3)}$ bundles on $S^6$ are classified by an
integer (and in fact the example above is a generator). However, using the
Eckmann-Hilton argument, one can show that over a sphere $S^{k+1}$, given a
$G$-bundle with transition function $t\colon S^k\to G$ representing a generator
of $g\in \pi_k(G)$, we can obtain the transition functions for the bundles
corresponding to elements $g^n$ by taking the pointwise power 
$t^n\colon S^k \to G$ for any $n\in \mathbb{Z}$. The same will be true for 
the lifted 2-bundles, where we take pointwise powers of the 2-group-valued
functor $(V^{[2]} \sourcetarget V) \to \String_H$. Thus, for spheres at least,
we can in principle give \v Cech cocycle descriptions for all String bundles.

As a final note, the abstract picture in the penultimate paragraph is
not restricted to smooth geometry: one can equally well take holomorphic
2-groups, assuming one has them. However, in current work with Raymond Vozzo
we have found that the basic gerbe on a simple, simply-connected complex
reductive Lie group, which is holomorphic \cite{Brylinski_94,Brylinski_00}, is also
multiplicative, so defines a weak 2-group in complex analytic stacks. This means
we can define holomorphic String bundles on complex homogeneous spaces, which
can be plugged into the higher twistor correspondence of Saemann-Wolf (eg \cite{Saemann-Wolf_14}).

\printbibliography

\end{document}